\documentclass[11pt]{amsart}

\newfont{\fnt}{cmr10 scaled 550}

\newtheorem{theorem}{Theorem}

\newtheorem{conj}{Conjecture}
\newtheorem{lemma}{Lemma}

\theoremstyle{remark}
\newtheorem{remark}{Remark}

\font\strange=msbm10

\renewcommand{\epsilon}{\varepsilon}

\renewcommand{\Sigma}{\varSigma}
\newcommand{\R}{{{\mathchoice  {\hbox{$\textstyle{\text{\strange R}}$}}
{\hbox{$\textstyle{\text{\strange R}}$}}
{\hbox{$\scriptstyle  N\kern-0.3em  R$}}  
{\hbox{$\scriptscriptstyle  R\kern-0.2em  R$}}}}}

\newcommand{\Z}{{{\mathchoice  {\hbox{$\textstyle{\text{\strange Z}}$}}
{\hbox{$\textstyle{\text{\strange Z}}$}}
{\hbox{$\scriptstyle  Z\kern-0.3em  Z$}}
{\hbox{$\scriptscriptstyle  Z\kern-0.2em  Z$}}}}}

\newcommand{\N}{{{\mathchoice  {\hbox{$\textstyle{\text{\strange N}}$}}
{\hbox{$\textstyle{\text{\strange N}}$}}
{\hbox{$\scriptstyle  N\kern-0.3em  N$}}
{\hbox{$\scriptscriptstyle  N\kern-0.2em  N$}}}}}

\renewcommand{\phi}{\varphi}
\usepackage{amsmath,amsthm,amscd}

\begin{document}
\title{On the DDVV Conjecture and the Comass in Calibrated Geometry (II)}
\date{August 20, 2007}

 \author{Zhiqin Lu}

\address{Department of
Mathematics, University of California,
Irvine, Irvine, CA 92697}


\email[Zhiqin Lu]{zlu@math.uci.edu}

\thanks{The
 author is partially supported by  NSF Career award DMS-0347033 and the
Alfred P. Sloan Research Fellowship.}

\maketitle 

\pagestyle{myheadings}

\newcommand{\ka}{K\"ahler }
\newcommand{\ii}{\sqrt{-1}}
\section{Introduction}
Let $M^n$ be an  $n$ dimensional manifold isometrically immersed into the space form $N^{n+m}(c)$ of constant sectional curvature $c$. Define the normalized scalar curvature $\rho$ and $\rho^\perp$ for the tangent bundle and the normal bundle as follows:
\begin{align}
\begin{split}
&\rho=\frac{2}{n(n-1)}\sum_{1=i<j}^n R(e_i, e_j,e_j,e_i),\\
&\rho^{\perp}=\frac{2}{n(n-1)}\left(\sum_{1=i<j}^n\sum_{1=r<s}^m\langle
R^\perp(e_i,e_j)\xi_r,\xi_s\rangle^2\right)^{\frac 12},
\end{split}
\end{align}
where $\{e_1,\cdots, e_n\}$ (resp. $\{\xi_1,\cdots,\xi_m\}$) is an orthonormal basis of the tangent (resp. normal space) at the point $x\in M$, and $R, R^\perp$ are the curvature tensors for the tangent and normal bundles, respectively.

In the study of submanifold theory, De Smet, Dillen, Verstraelen, and Vrancken~\cite{ddvv} made the following {\sl DDVV Conjecture}: 

\begin{conj}\label{cc1} Let $h$ be the second fundamental form, and 
let $H=\frac 1n \,{\rm trace}\, h$ be the mean curvature tensor.  Then 
\[
\rho+\rho^\perp\leq |H|^2+c.
\]
\end{conj}

 A weaker version of the above conjecture,
 \[
 \rho\leq |H|^2+c,
 \]
 was proved in ~\cite{ch1}. An alternate proof is in~\cite{bogd2}.
 
 In~\cite{ddvv}, the authors proved the following
 \begin{theorem}\label{ddvv}
 If $m=2$, then the conjecture is true.
 \end{theorem}
 
 In this  paper, we prove the conjecture in the case $m= 3$.    In the next version of this paper, we will prove $P(n,m)$.
 
 This paper is the continuation of the previous paper~\cite{cl}, where the case $n=3$ was proved.

 Let $x\in M$ be a fixed point and let $(h_{ij}^r)$ ($i,j=1,\cdots,n$ and $r=1,\cdots,m$) be the coefficients of the second fundamental form under some orthonormal basis. Then by
 Suceav{\u{a}}~\cite{bogd}, or ~\cite{dfv}, Conjecture~\ref{cc1} can be formulated as an inequality with respect to the coefficients $h_{ij}^r$ as follows:
  \begin{align}\label{1a}
\begin{split}
&\sum_{r=1}^m\sum_{1=i<j}^n(h_{ii}^r-h_{jj}^r)^2+2n\sum_{r=1}^m\sum_{1=i<j}^n(h_{ij}^r)^2\\
&\geq 2n\left(\sum_{1=r<s}^m\sum_{1=i<j}^n\left(\sum_{k=1}^n(h_{ik}^rh_{jk}^s-h_{ik}^sh_{jk}^r)\right)^2\right)^{\frac 12}.
\end{split}
\end{align}

Suppose that $A_1,A_2,\cdots,A_m$ are $n\times n$ symmetric real matrices. Let
\[
||A||^2=\sum_{i,j=1}^n a_{ij}^2,
\]
where $(a_{ij})$ are the entries of $A$, 
and let
\[
[A,B]=AB-BA
\]
be the commutator. Then the equation~\eqref{1a}, in terms of matrices, can be formulated as follows
\begin{conj}\label{conj2}
 For $n, m\geq 2$, we have
\begin{equation}\label{conj}
(\sum_{r=1}^m||A_r||^2)^2\geq 2(\sum_{r<s}||[A_r,A_s]||^2).
\end{equation}
Fixing $n,m$, we call the above inequality Conjecture $P(n,m)$.
\end{conj}

\begin{remark}\label{remk}
For derivation of ~\eqref{1a}, see~\cite[Theorem 2]{dfv}. Note that the prototype of the matrices are the traceless part of the second fundamental forms.
\end{remark}

{\bf Acknowledgment.} We thank Professor X-L Xin to bring us to the attention of the papers~\cite{lianmin,xusenlin}, where we learned on of the important techniques in this paper.

\section{Pinching theorems.}

Let $A_1,\cdots,A_m$ be $n\times n$ symmetric matrices. Let $P(n,m)$
 be the following conjecture~\cite{ddvv,cl}:
 
 \begin{conj} Using the above notations, we have
 \[
 2\sum_{i<j}||[A_i,A_j]||^2\leq \left(\sum_{i=1}^m||A_i||^2\right)^2.
 \]
 \end{conj}
 
 In ~\cite{lianmin,xusenlin}, the following result was proved (cf.~\cite[pp 585, equation (5)]{lianmin}):
 
 \begin{theorem} Using the same notations as above, we have
  \[
  2\sum_{i<j}||[A_i,A_j]||^2\leq \frac 32\left(\sum_{i=1}^m||A_i||^2\right)^2-\sum_{i=1}^m||A_i||^4.
  \]
  \end{theorem}
  
  \qed
  
  We denote  the above inequality to be $P'(n,m)$.  
  In this note, we prove the following
  
  \begin{theorem}\label{thm}
  \[
  P(n,m)\Rightarrow P'(n,m).
  \]
  \end{theorem}

  {\bf Proof.} We assume that 
  \[
  ||A_1||\geq\cdots \geq ||A_m||.
  \]
  We prove $P'(n,m)$ by induction: suppose $P'(n,m-1)$ is true. Then we have the following 
  
  \begin{lemma}
  If $P'(n,m)$ is true for 
  \[
  ||A_1||^2\leq\sum_{i=2}^m||A_i||^2,
  \]
  then $P'(n,m)$ is true for any $A_1,\cdots,A_m$.
  \end{lemma}
  
  {\bf Proof.} We let $A_1=tA_1'$ and assume that $||A'_1||=1$. Then inequality $P'(n,m)$ can be written as
  \begin{align}\label{b}
  \begin{split}&
  \frac 12 t^4-t^2\left (2\sum_{i=2}^m||[A_1',A_i]||^2-3\sum_{i=2}^m||A_i||^2\right)\\&
  +\frac 32\left(\sum_{i=2}^m||A_i||\right)^2-\sum_{i=2}^m||A_i||^4
  -2\sum_{2\leq i<j}||[A_i,A_j]||^2\geq 0.
  \end{split}
  \end{align}
 By the inductive assumption, the total of the last three terms of the above is nonnegative. Let
 \begin{equation}\label{defa}
  a=2\sum_{i=2}^m||[A_1',A_i]||^2-3\sum_{i=2}^m||A_i||^2.
  \end{equation}
  If $a\leq 0$, then
  then ~\eqref{b} is trivially true. On the other hand, if $a>0$, then the minimum value is obtained at
  \[
  t^2=a.
  \]
Using the fact that $||[A_1',A_i]||^2\leq 2||A_i||^2$, we obtain:
  \[
  ||A_1||^2\leq \sum_{i=2}^m||A_i||^2.
  \]
  
  \qed

  {\bf Proof of  Theorem~\ref{thm}.} If 
    \[
  ||A_1||^2\leq\sum_{i=2}^m||A_i||^2,
  \]
  then
  \[
  \left(\sum_{i=1}^m||A_i||^2\right)^2
  \leq   
  \frac 32\left(\sum_{i=1}^m||A_i||^2\right)^2-\sum_{i=1}^m||A_i||^4.
  \]
  Thus
  \[
 P(n,m)\Rightarrow P'(n,m).
   \]
  
Since $P(3,m)$ is true by the main result in~\cite{cl}, can we get new pinching constant using this new inequality?

\section{Proof of $P(n,3)$.}

In this section, we prove the following

\begin{theorem}\label{thm1}
Let $A,B,C$ be symmetric $n\times n$ matrices.  Then
\[
(||A||^2+||B||^2+||C||^2)^2\geq 2(||[A,B]||^2+||[B,C]||^2+||[C,A]||^2).
\]
\end{theorem}

We first prove the following lemma:

\begin{lemma} \label{lem1}  Let $x\geq y\geq 0$. Let $(\eta_1,\cdots,\eta_n)$ be  a unit vector. Then if $\{i,j\}\neq\{k,l\}$, we have
\[
(\eta_i-\eta_j)^2 x+(\eta_k-\eta_l)^2 y\leq 2x+y.
\]
\end{lemma}

{\bf Proof.} If $i\not\in\{k,l\}$ and $j\not\in\{k,l\}$, then we have
\[
(\eta_i-\eta_j)^2 x+(\eta_k-\eta_l)^2 y\leq 2(\eta_i^2+\eta_j^2) x+2(\eta_k^2+\eta_l^2)y
\leq 2(\eta_i^2+\eta_j^2) x+2(1-\eta_i^2-\eta_j^2) y.
\]
Thus we have
\[
(\eta_i-\eta_j)^2 x+(\eta_k-\eta_l)^2 y\leq 2(x-y)+2y=2x\leq 2x+y.
\]
On the other hand, if $i=k,l$ or $j=k,l$,
then WLOG, we can assume that $i=k=1,j=2,l=3$. Thus we have
\[
(\eta_1-\eta_2)^2x+(\eta_1-\eta_3)^2y=(\eta_1,\eta_2,\eta_3)
\begin{pmatrix}
x+y&-x&-y\\
-x&x&0\\
-y&0&y
\end{pmatrix}
\begin{pmatrix}
\eta_1\\\eta_2\\\eta_3
\end{pmatrix}.
\]
The largest eigenvalue of the above matrix is $x+y+\sqrt{x^2-xy+y^2}\leq 2x+y$. Since $\eta_1^2+\eta_2^2+\eta_3^2\leq 1$, we have
\[
(\eta_1-\eta_2)^2x+(\eta_1-\eta_3)^2y\leq 2x+y,
\]
as stated.

\qed

\begin{lemma}\label{lem2}
Suppose that $||A||^2+||B||^2+||C||^2=1$ and $||A||\geq||B||\geq||C||$. Let
\[
\lambda={\rm Max}\,(||[A,B]||^2+||[B,C]||^2+||[C,A]||^2),
\]
and let $A,B,C$ be the maximum point. Then we have
\[
2\lambda ||A||^2=||[A,B]||^2+||[A,C]||^2.
\]
\end{lemma}

{\bf Proof.} Consider the function
\[
F=||[A,B]||^2+||[B,C]||^2+||[C,A]||^2-\lambda'(||A||^2+||B||^2+||C||^2-1).
\]
Using the Lagrange multiplier's method, for any symmetric matrix $\xi$, we have
\begin{align*}
&\langle [A,B],[A,\xi]\rangle+\langle[C,B],[C,\xi]\rangle-\lambda'\langle B,\xi\rangle=0\\
&\langle [B,A],[B,\xi]\rangle+\langle[C,A],[C,\xi]\rangle-\lambda'\langle A,\xi\rangle=0\\
&\langle [B,C],[B,\xi]\rangle+\langle[A,C],[A,\xi]\rangle-\lambda'\langle C,\xi\rangle=0.
\end{align*}
Since $\xi$ is arbitrary, we have
\begin{align*}
&||[A,B]||^2+||[C,B]||^2-\lambda'||B||^2=0\\
&||[A,B]||^2+||[C,A]||^2-\lambda'||A||^2=0\\
&||[B,C]||^2+||[A,C]||^2-\lambda'||C||^2=0.
\end{align*}

Summing over the three equations, we have
\[
2\lambda=\lambda'.
\]
The lemma follows.

\qed

{\bf Proof of Theorem~\ref{thm1}.} Let 
\[
G=O(n)\times O(3).
\]
The group acts on $(A,B,C)$ as follows:
let $Q\in O(n)$, then the $Q$ action is
\[
(A,B,C)\mapsto (QAQ^T,QBQ^T,QCQ^T);
\]
let $Q_1=(q_{ij})\in O(3)$, then the $Q_1$ action is
\[
(A,B,C)\mapsto (q_{11}A+q_{12}B+q_{13}C,\cdots,q_{31}A+q_{32}B+q_{33}C).
\]
It is not hard to see that the inequality and the expression
\[
||[A,B]||^2+||[B,C]||^2+||[C,A]||^2
\]
are $G$ invariant. Thus WLOG, we assume that $A,B,C$ are orthogonal and consider the maximum of 
\[
||[A,B]||^2+||[A,C]||^2
\]
under the constraint $||B||^2=x, ||C||^2=y$ and $x\geq y$. We assume that $A$ is diagnolized.  Let $A'=A/||A||$, and let
\[
A'=\begin{pmatrix}
\eta_1\\
&\ddots\\
&& \eta_n
\end{pmatrix}
\]
Then $\eta_1^2+\cdots+\eta^2_n=1$.

Consider the function
\[
g=\sum_{i,j} (\eta_i-\eta_j)^2 (b_{ij}^2+c_{ij}^2)+\lambda_1(\sum_{i,j} b_{ij}^2-x)
+\lambda_2(\sum_{i,j}c_{ij}^2-y)+\mu(\sum_{i,j}b_{ij}c_{ij}).
\]

Using the Lagrange muliplier's method, at the maximum points, we have

\begin{align*}
&2((\eta_i-\eta_j)^2 +\lambda_1 )b_{ij}+\mu c_{ij}=0\\
&2((\eta_i-\eta_j)^2 +\lambda_2 )c_{ij}+\mu b_{ij}=0
\end{align*}
for $i\geq j$.

WLOG, we assume that $(\eta_i-\eta_j)^2$ are different. If $\mu=0$, then at most for one $i> j$ and one $k> l$, we have $b_{ij}\neq 0$ and $c_{kl}\neq 0$. Since $B,C$ are orthogonal, if $b_{ij}\neq 0$ and $c_{kl}\neq 0$, then we have
$(i,j)\neq (k,l)$. It follows that
\[
\sum_{i,j}(\eta_i-\eta_j)^2(b_{ij}^2+c_{ij}^2)\leq (\eta_i-\eta_j)^2x+(\eta_k-\eta_l)^2 y.
\]
By Lemma~\ref{lem1}, we have
\[
\sum_{i,j}(\eta_i-\eta_j)^2(b_{ij}^2+c_{ij}^2)\leq 2x+y.
\]
If $\mu\neq 0$, then we have
\[
((2(\eta_i-\eta_j)^2+\lambda_1)(2(\eta_i-\eta_j)^2+\lambda_2)-\mu^2)b_{ij}c_{ij}=0.
\]
Thus at most two pairs of $b_{ij}c_{ij}\neq 0$ for $i>j$. On the other hand, since $\mu\neq 0$, $b_{ij}\neq 0$ iff $c_{ij}\neq 0$. 
There are several cases:

{\bf Case 1.}
Suppose that $b_{ij}c_{ij}\neq 0$ and $b_{kl}c_{kl}\neq 0$ for $\{i,j\}\neq\{k,l\}$. Then $b_{ii}=c_{ii}=0$. The orthogonal condition implies that
\[
b_{ij}c_{ij}+b_{kl}c_{kl}=0.
\]
Using the above conditions, we can assume that
\[
b_{ij}=\sqrt {\frac x2}\cos\alpha, b_{kl}=\sqrt \frac x2\sin\alpha,c_{ij}=-\sqrt \frac y2\sin\alpha,c_{kl}=\sqrt \frac y2\cos\alpha.
\]
Thus 
\[
\sum_{i,j}(\eta_i-\eta_j)^2(b_{ij}^2+c_{ij}^2)=(\eta_i-\eta_j)^2(x\cos^2\alpha+y\sin^2\alpha)
+(\eta_k-\eta_l)^2(x\sin^2\alpha+y\cos^2\alpha).
\]
Apparently, the maximum values are obtained at $\alpha=0$ or $\frac\pi 2$. Using Lemma ~\ref{lem1}, in either case, we have
\[
\sum_{i,j}(\eta_i-\eta_j)^2(b_{ij}^2+c_{ij}^2)\leq 2x+y.
\]

{\bf Case 2.} If there is only one $b_{ij}\neq 0$, then we have
\begin{equation}\label{aac}
2b_{ij}c_{ij}+\sum_i b_{ii}c_{ii}=0.
\end{equation}
Thus we have
\[
\sum_{i,j}(\eta_i-\eta_j)^2(b_{ij}^2+c_{ij}^2)\leq 4(b_{ij}^2+c_{ij}^2).
\]
An element computation gives that
\[
4(b_{ij}^2+c_{ij}^2)\leq 2x+y
\]
using ~\eqref{aac}.

{\bf Case 3.} If $b_{ij}=0$ for $i>j$, then
\[
\sum_{i,j}(\eta_i-\eta_j)^2(b_{ij}^2+c_{ij}^2)=0\leq 2x+y.
\]

In summary, we have
\[
||[A',B]||^2+||[A',C]||^2\leq 2||B||^2+||C||^2.
\]

Using Lemma~\ref{lem2}, we have
\[
2\lambda||A||^2\leq ||A||^2 (2x+y).
\]
Since $||A||\geq||B||\geq ||C||$, we have $2x+y\leq 1$. Thus $2\lambda\leq 1$. This is what we want to prove.

\qed

\bibliographystyle{abbrv}   
\bibliography{pub,unp,2007}

\end{document}